\documentclass[11pt]{article}
\usepackage{amssymb}
\usepackage{latexsym,amsmath}
\newtheorem{precor}{{\bf Corollary}}

\newtheorem{precon}{{\bf Conjecture}}

\newenvironment{con}{\begin{precon}{\hspace{-0.5
               em}{\bf.\ }}}{\end{precon}}
\newtheorem{prealphcon}{{\bf Conjecture}}

\newtheorem{predefin}{{\bf Definition}}

\newtheorem{preexm}{{\bf Example}}

\newtheorem{preappl}{{\bf Application}}

\newtheorem{prelem}{{\bf Lemma}}

\newtheorem{preproof}{{\bf Proof.\ }}

\newenvironment{proof}[1]{\begin{preproof}{\rm
               #1}\hfill{$\blacksquare$}}{\end{preproof}}
\newtheorem{pretheorem}{{\bf Theorem}}

\newenvironment{theorem}{\begin{pretheorem}{\hspace{-0.5
               em}{\bf.\ }}}{\end{pretheorem}}
\newtheorem{prealphtheorem}{{\bf Theorem}}

\newtheorem{prealphlem}{{\bf Lemma}}

\newtheorem{prepro}{{\bf Proposition}}

\newenvironment{pro}{\begin{prepro}{\hspace{-0.5
               em}{\bf.\ }}}{\end{prepro}}
\newtheorem{preprb}{{\bf Problem}}

\newtheorem{prerem}{{\bf Remark}}

\newtheorem{preapp}{{\bf Application}}

\newtheorem{prequ}{{\bf Question}}

\newenvironment{qu}{\begin{prequ}{\hspace{-0.5
               em}{\bf.\ }}}{\end{prequ}}
\def\conct[#1,#2]{\mbox {${#1} \leftrightarrow {#2}$}}
\def\dconct[#1,#2]{\mbox {${#1} \rightarrow {#2}$}}
\def\deg[#1,#2]{\mbox {$d_{_{#1}}(#2)$}}
\def\mindeg[#1]{\mbox {$\delta_{_{#1}}$}}
\def\maxdeg[#1]{\mbox {$\Delta_{_{#1}}$}}
\def\outdeg[#1,#2]{\mbox {$d_{_{#1}}^{^+}(#2)$}}
\def\minoutdeg[#1]{\mbox {$\delta_{_{#1}}^{^+}$}}
\def\maxoutdeg[#1]{\mbox {$\Delta_{_{#1}}^{^+}$}}
\def\indeg[#1,#2]{\mbox {$d_{_{#1}}^{^-}(#2)$}}
\def\minindeg[#1]{\mbox {$\delta_{_{#1}}^{^-}$}}
\def\maxindeg[#1]{\mbox {$\Delta_{_{#1}}^{^-}$}}

\def\dre[#1,#2,#3]{\mbox {${\cal E}^{^{#3}}(#1,#2)$}}
\def\var[#1,#2]{\mbox {${\rm Var}_{_{#1}}(#2)$}}
\def\ls[#1]{\mbox {$\xi^{^{#1}}$}}
\def\hom[#1,#2]{\mbox {${\rm Hom}({#1},{#2})$}}
\def\onvhom[#1,#2]{\mbox {${\rm Hom^{v}}(#1,#2)$}}
\def\onehom[#1,#2]{\mbox {${\rm Hom^{e}}(#1,#2)$}}
\def\core[#1]{\mbox {$#1^{^{\bullet}}$}}
\def\cay[#1,#2]{\mbox {${\rm Cay}({#1},{#2})$}}
\def\sch[#1,#2,#3]{\mbox {${\rm Sch}({#1},{#2},{#3})$}}
\def\cays[#1,#2]{\mbox {${\rm Cay_{s}}({#1},{#2})$}}
\def\dirc[#1]{\mbox {$\stackrel{\rightarrow}{C}_{_{#1}}$}}
\def\cycl[#1]{\mbox {${\bf Z}_{_{#1}}$}}

\begin{document}

\begin{center} 
{\Large \bf The Matching Kneser Graph Conjecture For High Chromatic Numbers}\\
\vspace{0.3 cm}
{\bf Saeed Shaebani}\\
{\it School of Mathematics and Computer Science}\\
{\it Damghan University, Damghan, Iran}\\
{\tt shaebani@du.ac.ir}\\ \ \\
\end{center}
\begin{abstract}
\noindent In this paper, we show that
for any positive integers $r$, $k$, $\Theta$, and $\Gamma$ such that $k \geq 2$ and $r \geq k + \Gamma$,
there exists a connected graph $G$ for which
$$\begin{array}{llcr}
	\omega (G) = \chi (G) = k,  &   \chi \left( G , rK_2 \right) = \Theta , & {\rm and} & |E(G)| - {\rm ex}\left( G , rK_2 \right) = \Theta + \Gamma .
\end{array}$$




\end{abstract}
\section{Introduction}

The graphs considered in this paper are finite and simple.
Also, a matching of size $r$ is simply abbreviated by an {\it $r$-matching}, and it is denoted by $rK_2$.

Based on a graph $G$, and considering all $r$-matchings in $G$, Alishahi and Hajiabolhassan 
\cite{ALISHAHI2015186, ALISHAHI20172366, alishahi_hossein_2020} defined {\it the matching Kneser graph}, denoted by $\left( G , rK_2 \right)$, as a graph which includes each $r$-matching of $G$ as a node; and two $r$-matchings
$\left\{x_1 y_1 , x_2 y_2 , \dots , x_r y_r\right\}$
and
$\left\{z_1 w_1 , z_2 w_2 , \dots , z_r w_r\right\}$
are adjacent in $\left( G , rK_2 \right)$ whenever
$$\left\{x_1 y_1 , x_2 y_2 , \dots , x_r y_r\right\} \cap \left\{z_1 w_1 , z_2 w_2 , \dots , z_r w_r\right\} = \varnothing .$$
Indeed, each $r$-matching of $G$ is considered as a vertex of $\left( G , rK_2 \right)$;
and any two edge-disjoint $r$-matchings of $G$ establish an edge in $\left( G , rK_2 \right)$.

One of the most famous problems in graph theory goes back to determining the chromatic number of the matching Kneser graph $\left( lK_2 , rK_2 \right)$.
In 1955, Kneser \cite{kneser1955aufgabe} proved that if $l \geq 2r-1$, then
$\chi \left( lK_2 , rK_2 \right) \leq l-2r+2$; and conjectured that this inequality is indeed an eqiality; that is,
$\chi \left( lK_2 , rK_2 \right) = l-2r+2$.
This conjecture was settled by Lov{\'a}sz's in 1978 in a break-through paper
\cite{MR514625} which opened the doors of a new branch of Mathematics, called {\it Topological Combinatorics}.
One can refer to \cite{MR1988723, MR2057690} for more details about this new branch of Mathematics.

By an {\it $\left( rK_2 \right)$-free graph}, we mean a graph which includes no $rK_2$ as a subgraph.

For a graph $G$ and a positive integer $r$, {\it the generalized Turán number} ${\rm ex}\left( G , rK_2 \right)$ is defined as the greatest nonnegative integer $x$
for which an $\left( rK_2 \right)$-free spanning subgraph of $G$ with exactly $x$
edges exists. In fact, ${\rm ex}\left( G , rK_2 \right)$ is the maximum possible number of $|E(H)|$, where $H$ ranges over all $\left( rK_2 \right)$-free spanning subgraphs of $G$.

For a graph $G$ and a subset $A$ of $E(G)$, the symbol $G-A$ stands for a graph
with $V(G-A) := V(G)$ and $E(G-A) := E(G) \setminus A$.
Indeed, $G-A$ is the spanning subgraph of $G$ obtained by removing all edges in $A$.

Since ${\rm ex}\left( G , rK_2 \right)$ is the maximum number of edges over all $\left( rK_2 \right)$-free spanning subgraphs of $G$, one may easily deduce that
$$|E(G)| - {\rm ex}\left( G , rK_2 \right) = \min \left\{ |A| : A \subseteq E(G) \ {\rm and \ also} \ G-A \ {\rm is} \ \left( rK_2 \right) -{\rm free} \right\} .$$
In other words, $|E(G)| - {\rm ex}\left( G , rK_2 \right)$ is equal to the least number of edges in $E(G)$ whose deletion from $G$ results in a remaining spanning subgraph of $G$
with no $rK_2$.

Surprisingly \cite{ALISHAHI20172366, alishahi_hossein_2020}, the expression $|E(G)| - {\rm ex}\left( G , rK_2 \right)$
provides a sharp upper bound for $\chi \left( G , rK_2 \right)$; that is,
$$\chi \left( G , rK_2 \right) \leq |E(G)| - {\rm ex}\left( G , rK_2 \right).$$
More surprisingly, for many important classes of graphs, the equality
$\chi \left( G , rK_2 \right) = |E(G)| - {\rm ex}\left( G , rK_2 \right)$
holds \cite{ALISHAHI20172366, alishahi_hossein_2020}.
Alishahi and Hajiabbolhassan \cite{ALISHAHI20172366, alishahi_hossein_2020} provided
several interesting sufficient conditions whose occurrence imply this fantastic equality.
They \cite{alishahi_hossein_2020} also observed that for positive integers $l$ and $r$, we have
$$\chi \left( lK_2 , rK_2 \right) = l-2(r-1)  <  l-(r-1) = \left| E\left( lK_2  \right) \right| - {\rm ex}\left( lK_2 , rK_2 \right) $$
provided that $l \geq 2r-1$ and $r\geq 2$.
So, they found that the strict inequality
$\chi \left( G , rK_2 \right) < |E(G)| - {\rm ex}\left( G , rK_2 \right)$
could also happen. Since $lK_2$ is disconnected for $l \geq 2$, the equality
$\chi \left( G , rK_2 \right) = |E(G)| - {\rm ex}\left( G , rK_2 \right)$
fails for some disconnected graphs. They \cite{alishahi_hossein_2020}
conjectured that strictness may happen only for some disconnected graphs $G$;
or in other words, connectivity implies the equality.
\begin{con} (\cite{alishahi_hossein_2020}) \label{Conjecture1}
	For all connected graphs $G$, we have
	$$\chi \left( G , rK_2 \right) = |E(G)| - {\rm ex}\left( G , rK_2 \right) .$$
\end{con}
Obviously, this conjecture holds for all graphs $G$ with $|V(G)| < 2r$;
because in this case we have
$\chi \left( G , rK_2 \right) = |E(G)| - {\rm ex}\left( G , rK_2 \right) = 0 .$
So, each counterexample to Conjecture \ref{Conjecture1} must have an order
greater than or equal to $2r$.
Iradmusa \cite{Iradmusa2023} provided counterexamples to this conjecture whose
orders are equal to the best least possible $2r$.

\noindent A connected cubic graph is said to be a {\it snark} if it satisfies
the following two conditions simultaneously:
\begin{itemize}
	\item $\kappa ' (G) \geq 2$, where $\kappa ' (G)$ denotes the edge-connectivity of $G$.
	\item $\chi ' (G) =4$, where $\chi ' (G)$ is the chromatic-index of $G$.
\end{itemize}
Iradmusa \cite{Iradmusa2023} showed that for $r\geq 4$, all snarks $G$
of order $|V(G)|=2r$ satisfy
$\chi \left( G , rK_2 \right) = 1$ and $|E(G)| - {\rm ex}\left( G , rK_2 \right) = 3 $.
\begin{theorem} \cite{Iradmusa2023} \label{Iradmusa}
	If $r\geq 4$ and $G$ is a snark of order $2r$, then we have
	$$\chi \left( G , rK_2 \right) = 1 \ {\rm and} \ |E(G)| - {\rm ex}\left( G , rK_2 \right) = 3 .$$
\end{theorem}
It is worth pointing out that for $r \in \{1,2,3\}$, there is no snark of
order $2r$. So, the condition $r\geq 4$ in Theorem \ref{Iradmusa}
could also be replaced by $r \in \mathbb{N}$.

\noindent The first importance of Iradmusa's nice counterexample is its
order, which is the best least possible $2r$.
Also, the second importance of Iradmusa's counterexample is its regularity
of low degree 3.

\noindent In Iradmusa's counterexample, the following four statements hold for all
$r \geq 4$:
$$
\begin{array}[pos]{ll}
	\Lambda _1 := \chi \left( G , rK_2 \right) = 1, &  \Lambda _2 := |E(G)| - {\rm ex}\left( G , rK_2 \right) = 3 , \\
	       &     \\
	\Lambda _3 := \frac{|E(G)| - {\rm ex} \left(  G , rK_2  \right)}{\chi \left(  G , rK_2  \right)}  = 3 , & \Lambda _4 := \Bigl( |E(G)| - {\rm ex} \left(  G , rK_2  \right) \Bigr) -  \chi \left(  G , rK_2  \right) = 2 .
\end{array}
$$
Thus, the problem of whether $\Lambda _1 , \Lambda _2 , \Lambda _3 $, and $\Lambda _4 $
could attain other positive integers, would be of interest.
In this regard, some appropriate connected bipartite graphs $G$ were constructed
in \cite{SHAEBANI202487}.
\begin{theorem} \cite{SHAEBANI202487} \label{Saeed}
	If $r$, $\Theta$, and $\Gamma$ are arbitrary positive integers with $r\geq 3$ and
	$\Gamma \leq r-2$, then there exists a connected bipartite graph $G$ such that
	 $$\begin{array}{lcr}
	 	\chi \left( G , rK_2 \right) = \Theta & {\rm and} & |E(G)| - {\rm ex}\left( G , rK_2 \right) = \Theta + \Gamma .
	 \end{array}$$
\end{theorem}
\begin{theorem} \cite{SHAEBANI202487} \label{Saeed1}
	If $r$ and $\Theta$ are positive integers with $r \geq 3$, then there exists a
	tree $T$ of radius two which satisfies
	$$\begin{array}{lcr}
		\chi \left( T , rK_2 \right) = \Theta & {\rm and} & |E(T)| - {\rm ex}\left( T , rK_2 \right) = \Theta + r-2 .
	\end{array}$$
\end{theorem}
The following proposition is an immediate consequence of Theorem \ref{Saeed1}.
\begin{pro} \cite{SHAEBANI202487} \label{Saeed2}
	For any positive integer $\Theta$, there exists a sequence of trees
	$\left( T_r \right) _{r=3} ^{\infty}$
	such that $\chi \left(  T_r , rK_2  \right) = \Theta$ for all $r\geq 3$; and besides,
	$$\lim_{r \rightarrow \infty}  \Bigl( |E\left( T_r \right)| - {\rm ex} \left(  T_r , rK_2  \right) \Bigr) = +\infty .$$
\end{pro}

An anonymous referee of \cite{SHAEBANI202487} propounded the following interesting question.
\begin{qu} \cite{SHAEBANI202487} \label{newquestion}
	So far, known counterexamples to Conjecture \ref{Conjecture1} have chromatic numbers two (for bipartite graphs) and three (for snarks). What about greater integers? For any integer $k \geq 4$, does there exist a connected graph $G$ which is a counterexample to Conjecture \ref{Conjecture1} and satisfies $\chi (G) = k$?
\end{qu}

In this paper, we aim to answer Question \ref{newquestion}.

\section{The Main Result}

This section is devoted to provide an affirmative answer to Question \ref{newquestion}.
After the final version of \cite{SHAEBANI202487} was published online, the present author observed that by some light refinements of the Proof of Theorem \ref{Saeed} in \cite{SHAEBANI202487}, 
an appropriate answer to Question \ref{newquestion} may be achievable. Since \cite{SHAEBANI202487} was published online, it was not possible to add another new result in that paper \cite{SHAEBANI202487}. So, in this Section of the present paper, we are concerned with the new result; which is an answer to Question \ref{newquestion}.
The following theorem is the main result of this paper. Its proof is similar in spirit to the Proof of Theorem \ref{Saeed} in \cite{SHAEBANI202487}.

\begin{theorem}
	For any positive integers $r$, $k$, $\Theta$, and $\Gamma$ such that $k \geq 2$ and $r \geq k + \Gamma$,
	there exists a connected graph $G$ for which
	$$\begin{array}{llcr}
		\omega (G) = \chi (G) = k,  &   \chi \left( G , rK_2 \right) = \Theta , & {\rm and} & |E(G)| - {\rm ex}\left( G , rK_2 \right) = \Theta + \Gamma .
	\end{array}$$
\end{theorem}
\begin{proof}
	{
		Put $t:=(r-1)- \Gamma$.
		So, we have $k-1\leq t \leq r-2$.
		Also, put $$l:= \Theta + 2 \Gamma = \Theta + 2 (r-1-t) . $$
		Now, consider a connected graph $G$ with $V(G) := V_{1}  \cup  V_{2}$ where
		\begin{itemize}
			\item $V_1 := \{x_1 , x_2 , \dots , x_l\} \cup \left\{w_1 , w_2 , \dots , w_{t {l \choose {r-t}} + l   }     \right\} ,$
			\item $V_2 := \{y_1 , y_2 , \dots , y_l\} \cup \{z_1 , z_2 , \dots , z_t\} ;$
		\end{itemize}
		whose edge set $E(G)$ is defined $E(G) := E_1 (G) \cup E_2 (G) \cup E_3 (G)$ as follows:
		\begin{itemize}
			\item $E_1 (G) := \{x_1 y_1 , x_2 y_2 , \dots , x_l y_l\}$,
			\item $E_2 (G) := \{v z_j : v \in V_1 \ {\rm and} \ 1\leq j \leq t\}$,
			\item $E_3 (G) := \{z_i z_j : 1\leq i < j \leq k-1\}$.
		\end{itemize}
		Since each $r$-matching in $G$ has at least $r-t$ edges in $\{x_1 y_1 , x_2 y_2 , \dots , x_l y_l\}$, one could find a graph
		homomorphism from $\left( G , rK_2 \right)$ to $\left( lK_2 , (r-t)K_2 \right) $;
		which implies
		$$\chi \left( G , rK_2 \right)  \leq  \chi \left( lK_2 , (r-t)K_2 \right) .$$
		On the other hand, to each matching ${\cal M}$ of size $r-t$ in $\{x_1 y_1 , x_2 y_2 , \dots , x_l y_l\}$, we can assign an $r$-matching
		${\cal M} \cup \widetilde{{\cal M}}$ in $G$ with the following three properties:
		\begin{itemize}
			\item $\widetilde{{\cal M}}$ is a $t$-matching.
			\item Each edge of $\widetilde{{\cal M}}$ is incident with both of $\{z_1 , \dots , z_t\}$ and $\left\{w_1 , \dots , w_{t {l \choose {r-t}} + l   }     \right\}$.
			\item For any two distinct $(r-t)$-matchings ${\cal M}$ and ${\cal N}$ in
			$\{x_1 y_1 , x_2 y_2 , \dots , x_l y_l\}$, two $t$-matchings $\widetilde{{\cal M}}$ and
			$\widetilde{{\cal N}}$ are edge-disjoint.
		\end{itemize}
		Now, the assignment ${\cal M} \longmapsto {\cal M} \cup \widetilde{{\cal M}}$ defines a graph homomorphism
		from the graph
		$\left( lK_2 , (r-t)K_2 \right) $ to the graph $\left( G , rK_2 \right)$; which implies
		$$\chi \left( lK_2 , (r-t)K_2 \right) \leq \chi \left( G , rK_2 \right) .$$
		We conclude that
		$$\chi \left( G , rK_2 \right)  =  \chi \left( lK_2 , (r-t)K_2 \right) = \Theta .$$
		
		Suppose that one chooses $\Theta + \Gamma$ arbitrary edges from $\{x_1 y_1 , x_2 y_2 , \dots , x_l y_l\}$
		and then removes them from $G$. Doing so, in the resulting subgraph, $V_2$ has exactly $\Gamma + t$
		non-isolated vertices. Since $\Gamma + t = r-1$, the resulting subgraph contains no $rK_2$; and
		therefore,
		$$|E(G)| - {\rm ex}\left( G , rK_2 \right) \leq \Theta + \Gamma .$$
		We shall have established the theorem if we prove that
		$|E(G)| - {\rm ex}\left( G , rK_2 \right) \geq \Theta + \Gamma .$
		In this regard, it is sufficient to show that if $A$ is an arbitrary subset of $E(G)$ such that
		$|A| = \Theta + \Gamma - 1$, then the resulting subgraph of $G$ obtained by removing all
		edges of $A$ from $G$ still contains an $rK_2$.
		We denote the resulting subgraph by $G-A$.
		
		\noindent
		Since $\left|  \{x_1 y_1 , x_2 y_2 , \dots , x_l y_l\} \setminus A   \right| \geq \left|  \{x_1 y_1 , x_2 y_2 , \dots , x_l y_l\}  \right| - |A|
		= l - (\Theta + \Gamma - 1) = (\Theta + 2 \Gamma) - (\Theta + \Gamma - 1) = \Gamma + 1$,
		we find that $G-A$ contains at least $\Gamma + 1$ edges from $\{x_1 y_1 , x_2 y_2 , \dots , x_l y_l\}$.
		Without losing the generality, we may assume that $$\{x_1 y_1 , x_2 y_2 , \dots , x_{\Gamma + 1} y_{\Gamma + 1}\} \subseteq E(G-A).$$
		Since $|A| = \Theta + \Gamma - 1 \leq \Theta + 2 \Gamma = l$, we find that there exist $t$ vertices in
		$\left\{w_1 , w_2 , \dots , w_{t {l \choose {r-t}} + l   }     \right\}$
		that are incident with none of edges of $A$.
		Without losing the generality, we may assume that these mentioned vertices are
		$w_1 , w_2 , \dots , w_t$.
		So, $$\left\{  w_1 z_1 , w_2 z_2 , \dots , w_t z_t\right\} \subseteq E(G-A) .$$
		We conclude that 
		$$\{x_1 y_1 , x_2 y_2 , \dots , x_{\Gamma + 1} y_{\Gamma + 1}\} \cup \left\{  w_1 z_1 , w_2 z_2 , \dots , w_t z_t\right\} \subseteq E(G-A).$$
		Hence, $G-A$ contains a matching of size $\Gamma + 1 + t = r$; and we are done.
	}
\end{proof}

\bibliographystyle{plain}
\def\cprime{$'$} \def\cprime{$'$}

\end{document}